\documentclass[10pt,twocolumn,twoside]{article}
\usepackage[top=1in,bottom=1in,left=0.75in,right=0.75in,paperwidth=8.5in,paperheight=11in]{geometry}
\usepackage{fancyhdr}
\usepackage{amsmath}
\usepackage{graphicx}
\usepackage{newtxmath}
\usepackage{authblk}
\usepackage[square,longnamesfirst,numbers]{natbib}

\newcommand{\equationref}[1]{(\ref{#1})}

\newcommand{\secondderiv}[2]{\frac{\hdiff{#1}{2}}{\diff{#2}^2} - \frac{\diff{#1}}{\diff{#2}} \frac{\hdiff{#2}{2}}{\diff{#2}^2}}
\newcommand{\diff}[1]{\mathrm{d}#1}
\newcommand{\diffop}{\mathrm{d}}
\newcommand{\derivop}[2]{D_{#1}^{#2}}
\newcommand{\hdiff}[2]{\mathrm{d}^{#2}#1}
\newcommand{\dydx}{\frac{\dy}{\dx}}
\newcommand{\dxdy}{\frac{\dx}{\dy}}
\newcommand{\dx}{\diff{x}}
\newcommand{\dz}{\diff{z}}
\newcommand{\dq}{\diff{q}}
\newcommand{\dy}{\diff{y}}
\newcommand{\du}{\diff{u}}
\newcommand{\dt}{\diff{t}}

\title{Extending the Algebraic Manipulability of Differentials}
\author[1]{Jonathan Bartlett}

\author[2]{Asatur Zh. Khurshudyan}

\affil[1]{The Blyth Institute, jonathan.bartlett@blythinstitute.org}

\affil[2]{Institute of Mechanics, NAS of Armenia}

\date{\today}

\begin{document}
\maketitle

\begin{abstract}
Treating differentials as independent algebraic units have a long history of use and abuse.  
It is generally considered problematic to treat the derivative as a fraction of differentials rather than as a holistic unit acting as a limit, though for practical reasons it is often done for the first derivative.
However, using a revised notation for the second and higher derivatives will allow for the ability to treat differentials as independent units for a much larger number of cases.
\end{abstract}

\section{Introduction}

The calculus of variations has had a long, rich history, with many competing notations and interpretations.
The fluxion was the original concept of the derivative invented by Isaac Newton, and even had a notation similar to the modern Lagrange notation.
A competing notation for the derivative is the Leibniz notation, where the derivative is expressed as a ratio of differentials, representing arbitrarily small (possibly infinitesimal) differences in each variable.

The calculus was originally thought of as examining infinitely small quantities.
When these infinitely small quantities were put into ratio with each other, the result could potentially be within the reals (a likely result for smooth, continuous functions).
But, on their own, these infinitesimals were thought of as infinitely close to zero.

The concept of an infinitesimal caused a great deal of difficulty within mathematics, and therefore calculus was revised for the derivative to represent the limit of a ratio.
In such a conception, $\dx$ and $\dy$ are not really independent units, but, when placed in ratio with each other, represent the limit of that ratio as the changes get smaller and smaller.
However, many were not pleased with the limit notion, preferring to view $\dx$ and $\dy$ as distinct mathematical objects.

This question over the ontological status of differentials was somewhat paralleled by preferences in notation.
Those favoring the validity of infinitesimals generally preferred the Leibniz notation, where $\dx$ and $\dy$ are at least visually represented as individual units, while those favoring the limit conception of the derivative generally prefer the Lagrange notation, where the derivative is a holistic unit.

In an interesting turn of events, in the late 19th century, the Leibniz \emph{notation} for the derivative largely won out, but the Langrangian \emph{conception} of the derivative has been the favored intellectual interpretation of it.
Essentially, this means that equations are generally written as if there were distinct differentials available, but they are manipulated as if they only represent limits of a ratio which cannot be taken apart.

This dichotomy has led to an unfortunate lack of development of the notation.
Because it is generally assumed that differentials are not independent algebraic units, the fact that issues arise when treating them as such has not caused great concern, and has simply reinforced the idea that they should not be treated algebraically.
Therefore, there has been little effort to improve the notation to allow for a more algebraic treatment of individual differentials.

However, as will be shown, the algebraic manipulability of differentials can be greatly expanded if the notation for higher-order derivatives is revised.
This leads to an overall simplification in working with calculus for both students and practitioners, as it allows items which are written as fractions to be treated as fractions.
It prevents students from making mistakes, since their natural inclination is to treat differentials as fractions.\footnote{Since many in the engineering disciplines are not formally trained mathematicians, this also can prevent professionals in applied fields from making similar mistakes.}
Additionally, there are several little-known but extremely helpful formulas which are straightforwardly deducible from this new notation.

Even absent these practical concerns, we find that reconceptualizing differentials in terms of algebraically-manipulable terms is an interesting project in its own right, and perhaps may help us see the derivative in a new way, and adapt it to new uses in the future.
There may also be additional formulas which can in the future be more directly connected to the algebraic formulation of the derivative.

\section{The Problem of Manipulating Differentials Algebraically}
\label{secproblem}

When dealing with the first derivative, there are generally few practical problems in treating differentials algebraically.
If $y$ is a function of $x$, then $\dydx$ is the first derivative of $y$ with respect to $x$.
This can generally be treated as a fraction.

For instance, since $\dxdy$ is the first derivative of $x$ with respect to $y$, it is easy to see that these values are merely the inverse of each other.
The inverse function theorem of calculus states that $\dxdy = \frac{1}{\dydx}$.
The generalization of this theorem into the multivariable domain essentially provides for fraction-like behavior within the first derivative.

Likewise, in preparation for integration, both sides of the equation can be multiplied by $\dx$.
Even in multivariate equations, differentials can essentially be multiplied and divided freely, as long as the manipulations are dealing with the first derivative.

Even the chain rule goes along with this.
Let $x$ depend on parameter $u$.  If one has the derivative $\frac{\dy}{\du}$ and multiplies it by the derivative $\frac{\du}{\dx}$ then the result will be $\frac{\dy}{\dx}$.
This is identical to the chain rule in Lagrangian notation.

It is well recognized that problems occur when if one tries to extend this technique to the second derivative \citep{swokowski1983}.
Take for a simple example the function $y = x^3$.
The first derivative is $\dydx = 3x^2$.
The second derivative is $\frac{\hdiff{y}{2}}{\dx^2} = 6x$.

Say that it is later discovered that $x$ is a function of $t$ so that $x = t^2$.
The problem here is that the chain rule for the second derivative is not the same as what would be implied by the algebraic representation.

Here we arrive at one of the major problematic points for using the current notation of the second derivative algebraically.
To demonstrate the problem explicitly, if one were to take the second derivative seriously as a set of algebraic units, one should be able to multiply $\frac{\hdiff{y}{2}}{\dx^2}$ by $\frac{\dx^2}{\dt^2}$ to get the second derivative of $y$ with respect to $t$.
However, this does not work.
If the differentials are being treated as algebraic units, then $\frac{\dx^2}{\dt^2}$ is the same as $\left(\frac{\dx}{\dt}\right)^2$, which is just the first derivative of $x$ with respect to $t$ squared.
The first derivative of $x$ with respect to $t$ is $\frac{\dx}{\dt} = 2t$.
Therefore, treating the second derivative algebraically would imply that all that is needed to do to convert the second derivative of $y$ with respect to $x$ into the second derivative of $y$ with respect to $t$ is to multiply by $(2t)^2$.

However, this reasoning leads to the false conclusion that $\frac{\hdiff{y}{2}}{\dt^2} = 24t^4$.
If, instead, the substitution is done at the beginning, it can be easily seen that the result should be $30t^4$:
\begin{align*}
y &= x^3 \\
x &= t^2 \\
y &= (t^2)^3 \\
y &= t^6 \\
y' &= 6t^5 \\
y'' &= 30t^4 
\end{align*}

This is also shown by the true chain rule for the second derivative, based on Fa\`a di Bruno's formula \citep{johnson2002}.
This formula says that the chain rule for the second derivative should be:
\begin{equation}
\label{chainrulesecond}
\frac{\hdiff{y}{2}}{\dt^2} = \frac{\hdiff{y}{2}}{\dx^2}\left(\frac{\dx}{\dt}\right)^2 + \dydx \frac{\hdiff{x}{2}}{\dt^2} 
\end{equation}
This, however, is extremely unintuitive, and essentially makes a mockery out of the concept of using the differential as an algebraic unit.

It is generally assumed that this is a problem for the idea that second differentials should be treated as algebraic units.
However, it is possible that the real problem is that the notation for second differentials has not been given as careful attention as it should.

The habits of mind that have come from this have even affected nonstandard analysis, where, despite their appreciation for the algebraic properties of differentials, have left the algebraic nature of the second derivative either unexamined (as in \citep{henle1979})  or examined poorly (i.e., leaving out the problematic nature of the second derivative, as in  \citep[pg.~4]{keisler2000}).

\section{A Few Notes on Differential Notation}
\label{secdiffnotation}

Most calculus students glaze over the notation for higher derivatives, and few if any books bother to give any reasons behind what the notation means.
It is important to go back and consider why the notation is what it is, and what the pieces are supposed to represent.

In modern calculus, the derivative is always taken with respect to some variable.
However, this is not strictly required, as the differential operation can be used in a context-free manner.
The processes of taking a \emph{differential} and solving for a \emph{derivative} (i.e., some ratio of differentials) can be separated out into logically separate operations.\footnote{The 
  idea that finding a differential (i.e., similar to a derivative, but not being with respect to any particular variable) 
  can be separated from the operation of finding a derivative (i.e., differentiating \emph{with respect to some particular variable}) 
  is considered an anathema to some, but this concept can be inferred directly from the activity of treating derivatives as fractions of differentials.  
  The rules for taking a differential are identical to those for taking an implicit derivative, but simply leaving out dividing the final differential by the differential of the independent variable.

  For those uncomfortable with taking a differential without a derivative (i.e., without specifying an independent variable), imagine the differential operator $\diffop()$ as combining the operations of taking an implicit derivative with respect to a non-present variable (such as $q$) followed by a multiplication by the differential of that variable (i.e., $\dq$ in this example).  
  So, taking the differential of $e^x$ is written as $\diffop(e^x)$ and the result of this operation is $e^x\,\dx$.  
  This is the same as if we had taken the derivative with respect to the non-present variable $q$ and then multiplied by $\dq$.  
  So, for instance, taking the differential of the function $e^x$, the operation would start out with a derivative with respect to $q$ $\frac{\diffop}{\dq}(e^x) = e^x \frac{\dx}{\dq}$ followed by a multiplication by $\dq$, yielding just $e^x\,\dx$.  

  Doing this yields the standard set of differential rules, but allows them to be applied separately from (and prior to) a full derivative.  
  Also note that because they have no dependency on any variable present in the equation, the rules work in the single-variable and multi-variable case.  
  Solving for a derivative is then merely solving for a ratio of differentials that arise after performing the differential.  
  It unifies explicit and implicit differentiation into a unified process that is easier to teach, use, and understand, and requires few if any special cases, save the standard requirements of continuity and smoothness.
}

In such an operation, instead of doing $\frac{\diffop}{\dx}$ (taking the derivative with respect to the variable $x$), one would separate out performing the differential and dividing by $\dx$ as separate steps.
Originally, in the Leibnizian conception of the differential, one did not even bother solving for derivatives, as they made little sense from the original geometric construction of them \citep[pgs.~8,~59]{bos1974}.

For a simple example, the differential of $x^3$ can be found using a basic differential operator such that $\diffop(x^3) = 3x^2\,\dx$.
The derivative is simply the differential divided by $\dx$.  
This would yield $\frac{\diffop(x^3)}{\dx} = 3x^2$.

For implicit derivatives, separating out taking the differential and finding a particular derivative greatly simplifies the process.
Given an function (say, $z^2 = \sin(q)$), the differential can be applied to both sides just like any other algebraic manipulation:
\begin{align*}
z^2 &= \sin(q) \\
\diffop(z^2) &= \diffop(\sin(q)) \\
2z\,\dz &= \cos(q)\,\dq
\end{align*}
From there, the equation can be manipulated to solve for $\frac{\dz}{\dq}$ or $\frac{\dq}{\dz}$, or it can just be left as-is.

The basic differential of a variable is normally written simply as $\diffop(x) = \dx$.
In fact, $\dx$ can be viewed merely a shorthand for $\diffop(x)$.

The second differential is merely the differential operator applied twice \citep[pg.~17]{bos1974}:
\begin{equation}
\diffop(\diffop(x)) = \diffop(\dx) = \hdiff{x}{2}
\end{equation}
Therefore, the second differential of a function is merely the differential operator applied twice.
However, one must be careful when doing this, as the product rule affects products of differentials as well.

For instance, $\diffop(3x^2\,\dx)$ will be found using the product rule, where $u = 3x^2$ and $v = \dx$.
In other words:
\begin{align*}
\diffop(3x^2\,\dx) &= 3x^2(\diffop(\dx)) + \diffop(3x^2)\,\dx \\
                 &= 3x^2\,\hdiff{x}{2} + 6x\,\dx\,\dx \\
                 &= 3x^2\,\hdiff{x}{2} + 6x\,\dx^2
\end{align*}
The point of all of this is to realize that the notation $\frac{\hdiff{y}{2}}{\dx^2}$ is not some arbitrary arrangement of symbols, but has a deep (if, as will be shown, slightly incorrect or misleading) meaning.
The notation means that the equation is showing the ratio of the second differential of $y$ (i.e., $\diffop(\diffop(y))$) to the square of $\dx$ (i.e., $\dx^2$).\footnote{
In Leibniz notation, $\dx^2$ is equivalent to $(\dx)^2$.  If the differential of $x^2$ was wanted, it would be written as $\diffop(x^2)$.  The rules are given in \citep[pg.~24]{bos1974}.}

In other words, starting with $y$, then applying the differential operator twice, and then dividing by $\dx$ twice, arrives at the result $\frac{\hdiff{y}{2}}{\dx^2}$.
Unfortunately, that is not the same sequence of steps that happens when two derivatives are performed, and thus it leads to a faulty formulation of the second derivative.

\section{Extending the Second Derivative's Algebraic Manipulability}
\label{secextendingmanipulability}

As a matter of fact, order of operations is very important when doing derivatives.
When doing a derivative, one \emph{first} takes the differential and \emph{then} divides by $\dx$.
The second derivative is the derivative of the first, so the next differential occurs \emph{after the first derivative is complete}, and the process finishes by dividing by $\dx$ again.

However, what does it look like to take the differential of the first derivative?
Basic calculus rules tell us that the quotient rule should be used:
\begin{align*}
\diff\left(\dydx\right) &= \frac{\dx(\diffop(\dy)) - \dy(\diffop(\dx))}{(\dx)^2} \\
                        &= \frac{\dx\,\hdiff{y}{2} - \dy\,\hdiff{x}{2}}{\dx^2} \\
                        &= \frac{\dx\,\hdiff{y}{2}}{\dx^2} - \frac{\dy\,\hdiff{x}{2}}{\dx^2} \\
                        &= \frac{\dx}{\dx}\frac{\hdiff{y}{2}}{\dx} - \frac{\dy}{\dx}\frac{\hdiff{x}{2}}{\dx} \\
                        &= \frac{\hdiff{y}{2}}{\dx} - \frac{\dy}{\dx}\frac{\hdiff{x}{2}}{\dx}
\end{align*}
Then, for the second step, this can be divided by $\dx$, yielding:
\begin{equation}
\label{secondderiv}
\frac{\diffop\left(\frac{\dy}{\dx}\right)}{\dx} = \frac{\hdiff{y}{2}}{\dx^2} - \frac{\dy}{\dx}\frac{\hdiff{x}{2}}{\dx^2}
\end{equation}
This, in fact, yields a notation for the second derivative which is equally algebraically manipulable as the first derivative.
It is not very pretty or compact, but it works algebraically.

The chain rule for the second derivative fits this algebraic notation correctly, provided we replace each instance of the second derivative with its full form (cf.~\equationref{chainrulesecond}):
\begin{equation}
\label{chainrulesecondiff} 
\resizebox{0.9\columnwidth}{!}{%
$\frac{\hdiff{y}{2}}{\diff{t}^2} - \frac{\diff{y}}{\diff{t}} \frac{\hdiff{x}{2}}{\diff{x}^2} = \left(\secondderiv{y}{x}\right)\left(\frac{\dx}{\dt}\right)^2 + \dydx \left(\secondderiv{x}{t}\right)$
}
\end{equation}
This in fact works out perfectly algebraically.

One objection that has been given to the present authors by early reviewers about the formula presented in \equationref{secondderiv} is that the ratio $\frac{\hdiff{x}{2}}{\dx^2}$ reduces to zero.
However, this is not necessarily true.
The concern is that, since $\frac{\dx}{\dx}$ is always $1$ (i.e., a constant), then $\frac{\hdiff{x}{2}}{\dx^2}$ should be zero.
The problem with this concern is that we are no longer taking $\frac{\hdiff{x}{2}}{\dx^2}$ to be the derivative of $\frac{\dx}{\dx}$.
Using the notation in \equationref{secondderiv}, the derivative of $\frac{\dx}{\dx}$ would be:
\begin{equation}
\label{dxdxderiv}
\frac{\diffop\left(\frac{\dx}{\dx}\right)}{\dx} = \frac{\hdiff{x}{2}}{\dx^2} - \frac{\dx}{\dx}\frac{\hdiff{x}{2}}{\dx^2}
\end{equation}
In this case, since $\frac{\dx}{\dx}$ reduces to $1$, the expression is obviously zero.
However, in \equationref{dxdxderiv}, the term $\frac{\hdiff{x}{2}}{\dx^2}$ is not itself necessarily zero, since it is \emph{not} the second derivative of $x$ with respect to $x$.

\section{The Notation for the Higher Order Derivatives}

The notation for the third and higher derivatives can be found using the same techniques as for the second derivative.
To find the third derivative of $y$ with respect to $x$, one starts with the second derivative and takes the differential:
\begin{align*}
\diffop&\left(\frac{\diffop\left(\dydx\right)}{\dx}\right) \\
&= \diffop\left(\secondderiv{y}{x}\right) \\
  &= \diffop\left(\frac{\dx\,\hdiff{y}{2} - \dy\,\hdiff{x}{2}}{\dx^3}\right) \\
  &= \frac{(\dx^3)(\diffop(\dx\,\hdiff{y}{2} - \dy\,\hdiff{x}{2})) - (\dx\,\hdiff{y}{2} - \dy\,\hdiff{x}{2})(\diffop(\dx^3))}{(\dx^3)^2} \\
%  &= \frac{(\dx^3)(\dx\,\diffop(\hdiff{y}{2}) + \hdiff{y}{2}\,\diffop(\dx) - (\dy\,\diffop(\hdiff{x}{2}) + \hdiff{x}{2}\,\diffop(\dy))) - (\dx\,\hdiff{y}{2} - \dy\,\hdiff{x}{2})(3\,\dx^2\,\hdiff{x}{2})}{\dx^6} \\
%  &= \frac{(\dx^3)(\dx\,\hdiff{y}{3} + \hdiff{y}{2}\,\hdiff{x}{2} - \dy\,\hdiff{x}{3} - \hdiff{x}{2}\,\hdiff{y}{2}) - (\dx\,\hdiff{y}{2} - \dy\,\hdiff{x}{2})(3\,\dx^2\,\hdiff{x}{2})}{\dx^6} \\
%  &= \frac{\dx^4\,\hdiff{y}{3} + \dx^3\,\hdiff{y}{2}\,\hdiff{x}{2} - \dx^3\,\dy\,\hdiff{x}{3} - \dx^3\,\hdiff{x}{2}\,\hdiff{y}{2} - 3\,\dx^3\,\hdiff{x}{2}\,\hdiff{y}{2} + 3\,\dx^2\,\dy\,(\hdiff{x}{2})^2}{\dx^6} \\
%  &= \frac{\dx^4\,\hdiff{y}{3} - \dx^3\,\dy\,\hdiff{x}{3} - 3\,\dx^3\,\hdiff{x}{2}\,\hdiff{y}{2} + 3\,\dx^2\,\dy\,(\hdiff{x}{2})^2}{\dx^6} \\
%  &= \frac{\dx^4\,\hdiff{y}{3}}{\dx^6} - \frac{\dx^3\,\dy\,\hdiff{x}{3}}{\dx^6} - \frac{3\,\dx^3\,\hdiff{x}{2}\,\hdiff{y}{2}}{\dx^6} + \frac{3\,\dx^2\,\dy\,(\hdiff{x}{2})^2}{\dx^6} \\
  &= \frac{\hdiff{y}{3}}{\dx^2} - \frac{\dy}{\dx}\,\frac{\hdiff{x}{3}}{\dx^2} - 3\frac{\hdiff{x}{2}}{\dx^2}\frac{\hdiff{y}{2}}{\dx} + 3\frac{\dy}{\dx}\,\frac{(\hdiff{x}{2})^2}{\dx^3} 
\end{align*}
Finally, this result is divided by $\dx$:
\begin{equation}
\resizebox{0.9\columnwidth}{!}{%
$\frac{\diffop\left(\frac{\diffop\left(\dydx\right)}{\dx}\right)}{\dx} = \frac{\hdiff{y}{3}}{\dx^3} - \frac{\dy}{\dx}\,\frac{\hdiff{x}{3}}{\dx^3} - 3\frac{\hdiff{x}{2}}{\dx^2}\frac{\hdiff{y}{2}}{\dx^2} + 3\frac{\dy}{\dx}\,\frac{(\hdiff{x}{2})^2}{\dx^4}$
}
\end{equation}

This expression includes a lot of terms not normally seen, so some explanation is worthwhile.
In this expression, $\hdiff{x}{2}$ represents the second differential of $x$, or $\diffop(\diffop(x))$. 
Therefore, $(\hdiff{x}{2})^2$ represents $(\diffop(\diffop(x)))^2$.
Likewise, $\dx^4$ represents $(\diffop(x))^4$.

Because the expanded notation for the second and higher derivatives is much more verbose than the first derivative, it is often useful to adopt a slight modification of Arbogast's $D$ notation (see \citep[pgs.~209,218--219]{cajori2}) for the total derivative instead of writing it as algebraic differentials:\footnote{The difference between this notation and that of Arbogast is that we are subscripting the $D$ with the variable with which the derivative is being taken with respect to. Additionally, we are always supplying in the superscript the number of derivatives we are taking.  Therefore, where Arbogast would write simply $D$, this notation would be written as $D_x^1$.}
\begin{align}
D_{x}^2y &= \secondderiv{y}{x} \\
D_{x}^3y &= \frac{\hdiff{y}{3}}{\dx^3} - \frac{\dy}{\dx}\,\frac{\hdiff{x}{3}}{\dx^3} - 3\frac{\hdiff{x}{2}}{\dx^2}\frac{\hdiff{y}{2}}{\dx^2} + 3\frac{\dy}{\dx}\,\frac{(\hdiff{x}{2})^2}{\dx^4}
\end{align}
This gets even more important as the number of derivatives increases.
Each one is more unwieldy than the previous one.
However, each level can be interconverted into differential notation as follows:
\begin{equation}
D_{x}^{n}y = \frac{\diffop(D_{x}^{n-1}y)}{\dx}
\end{equation}
The advantages of Arbogast's notation over Lagrangian notation are that (1) this modification of Arbogast's notation clearly specifies both the top and bottom differential, and (2) for very high order derivatives, Lagrangian notation takes up $n$ superscript spaces to write for the $n$th derivative, while Arbogast's notation only takes up $\log(n)$ spaces.

Therefore, when a compact representation of higher order derivatives is needed, this paper will use Arbogast's notation for its clarity and succinctness.\footnote{It may be surprising to find a paper on the algebraic notation of differentials using a non-algebraic notation.  The goal, however, is to only use ratios \emph{when they act as ratios}.  When writing a ratio that works like a ratio is too cumbersome, we prefer simply avoiding the ratio notation altogether, to prevent making unwarranted leaps based on notation that may mislead the intuition.}

\section{Swapping the Independent and Dependent Variables}
\label{secswapping}

In fact, just as the algebraic manipulation of the first derivative can be used to convert the derivative of $y$ with respect to $x$ into the derivative of $x$ with respect to $y$, combining it with Arbogast's notation for the second derivative can be used to generate the formula for doing this on the second derivative:
\begin{align*}
\derivop{x}{2}y &= \secondderiv{y}{x} \\
\derivop{x}{2}y \frac{\dx^3}{\dy^3} &= \frac{\hdiff{y}{2}}{\dx^2} \frac{\dx^3}{\dy^3} - \frac{\dy}{\dx}\frac{\hdiff{x}{2}}{\dx^2} \frac{\dx^3}{\dy^3} \\
\derivop{x}{2}y \left(\frac{\dx}{\dy}\right)^3 &= \frac{\hdiff{y}{2}}{\dy^2} \frac{\dx}{\dy} - \frac{\hdiff{x}{2}}{\dy^2} \\
- \derivop{x}{2}y \left(\frac{\dx}{\dy}\right)^3 &= \frac{\hdiff{x}{2}}{\dy^2} - \frac{\dx}{\dy} \frac{\hdiff{y}{2}}{\dy^2} \\
- \derivop{x}{2}y \left(\frac{1}{\frac{\dy}{\dx}}\right)^3 &= \frac{\hdiff{x}{2}}{\dy^2} - \frac{\dx}{\dy}\frac{\hdiff{y}{2}}{\dy^2}  \\
- \derivop{x}{2}y \left(\frac{1}{\derivop{x}{1}y}\right)^3 &= \frac{\hdiff{x}{2}}{\dy^2} - \frac{\dx}{\dy}\frac{\hdiff{y}{2}}{\dy^2} 
\end{align*}
It can be seen that this final equation is the derivative of $x$ with respect to $y$.
Therefore, it can generally be stated that the second derivative of $y$ with respect to $x$ can be transformed into the second derivative of $x$ with respect to $y$ with the following formula:
\begin{equation}\label{inverse}
- \derivop{x}{2}y \left(\frac{1}{\derivop{x}{1}y}\right)^3 = \derivop{y}{2}x
\end{equation}
To see this formula in action on a simple equation, consider $y = x^3$.  
Performing two derivatives gives us:
\begin{align}
y &= x^3 \label{inverseexorig} \\
\derivop{x}{1}y &= 3x^2 \\
\derivop{x}{2}y &= 6x
\end{align}
According to \equationref{inverse}, $\derivop{y}{2}x$ (or, $x''$ in Lagrangian notation) can be found by performing the following:
\begin{align}
\derivop{y}{2}x &= -(6x) \left(\frac{1}{3x^2}\right)^3 \nonumber \\
 &= \frac{-6x}{27x^6} \nonumber \\
 &= \frac{-2}{9}x^{-5} \label{inverseexsecondbyformula}
\end{align}
This can be checked by taking successive derivatives of the inverse function of \equationref{inverseexorig}:
\begin{align}
x &= y^{\frac{1}{3}} \nonumber \\
\derivop{y}{1}x &= \frac{1}{3}y^{\frac{-2}{3}} \nonumber \\
\derivop{y}{2}x &= -\frac{2}{9}y^{\frac{-5}{3}} \label{inverseexsecond}
\end{align}
\equationref{inverseexsecond} can be seen to be equivalent to \equationref{inverseexsecondbyformula} by substituting for $y$ using \equationref{inverseexorig}:
\begin{align}
\derivop{y}{2}x &= -\frac{2}{9}(x^3)^{\frac{-5}{3}} \nonumber \\
 &= -\frac{2}{9}x^{-5} 
\end{align}
This is the same result achieved by using the inversion formula (cf. \equationref{inverse}).

%A similar process can be done to find the formula for swapping dependent and independent variables with the third derivative:
%
%\begin{align*}
%\derivop{x}{3}y &= \frac{\hdiff{y}{3}}{\dx^3} - \frac{\dy}{\dx}\,\frac{\hdiff{x}{3}}{\dx^3} - 3\frac{\hdiff{x}{2}}{\dx^2}\frac{\hdiff{y}{2}}{\dx^2} - 3\frac{\dy}{\dx}\,\frac{(\hdiff{x}{2})^2}{\dx^4}
%\end{align*}
%FIXME

\section{Using the Inversion Formula for the Second Derivative}
\label{secinversionapplication}

While the inversion formula (cf. \equationref{inverse}) is not original, it is a tool that many mathematicians are unaware of, and is rarely considered for solving higher-order differentials.\footnote{The authors of this paper, as well as several early reviewers, had originally thought that the inversion formula was a new finding.
  Again, that is the usefulness of the notation.
  Specific formulas such as the inversion formula do not need to be taught, as they simply flow naturally out of the notation.
  Even though the inversion formula is not new with this paper, showing how the present authors were able to use it to good benefit demonstrates the benefit of an improved notation---practitioners needs not memorize endless formulas, but they can be developed straightforwardly as needed based upon basic intuitions.
}

As an example of how to apply \equationref{inverse}, consider
second order ordinary nonlinear differential equations of the form
\[
F\left(y'', y', y\right) = 0.
\]
Equations of this form can be solved implicitly for
\[
F\left(a, b, c\right) = a - b^3 f(c)
\]
for generic function $f$. 
Indeed, consider the equation
\begin{equation}\label{ODE}
\derivop{x}{2} y = f(y) \left(\frac{\dy}{\dx}\right)^3.
\end{equation}
Then, by virtue of \equationref{inverse} we derive
\[
\derivop{y}{2} x = - f(y).
\]
Integration of this equation with respect to $y$ twice will provide with
\begin{equation}\label{implicitsol}
x(y) = - \int \int f(y) \, \dy\, \dy.
\end{equation}

For simplicity, let
\[
f(y) = y,
\]
so that \equationref{ODE} is reduced to
\[
\derivop{x}{2} y = y \left(\frac{d y}{dx}\right)^3,
\]
the real exact solutions of which is
\begin{equation}\label{sol}\begin{split}
y(x) &= \frac{6 \sqrt[3]{2} c_1}{\sqrt[3]{162(x + c_2) + \sqrt{23328 c_1^3 + \left[162 (x + c_2) \right]^2}}} - \\
&- \frac{\sqrt[3]{162(x + c_2) + \sqrt{23328 c_1^3 + \left[162 (x + c_2) \right]^2}}}{3 \sqrt[3]{2}}
\end{split}\end{equation}
Here $c_1$ and $c_2$ are integration constants that must be determined from given boundary or Cauchy conditions.

On the other hand, \equationref{implicitsol} results in
\[
x(y) = \frac{y^3}{6} + c_1 y + c_2,
\]
the real inverse of which exactly coincides with \equationref{sol}.

\section{Relationship to Historic Leibnizian Thought}
\label{secliebnizhist}

The view of differentials presented by Leibniz and those following in his footsteps differed significantly from the modern-day view of calculus.
The modern view of calculus focuses on functions, which have defined independent and dependent variables.
The Leibniz view, however, according to \citep{bos1974}, is a much more geometric view.
There is no preferred independent or dependent variable.

The modern concept of the derivative generally implies a dependent and in independent variable.
The numerator is the dependent variable and the denominator is the independent variable.
In the geometric view, however, there are only relationships, and these relationships do not necessarily have an implied dependency relationship.

Therefore, Leibnizian differentiation doesn't occur with respect to any independent variable.
There is no preferred independent variable.
Likewise, as we have seen in Sections~\ref{secswapping} and~\ref{secinversionapplication}, the version of the differential presented here allows for the reversal of variable dependency relationships.
Similarly, the procedure of differentiation given in Section~\ref{secdiffnotation} which allows us to formulate the new notation for the second derivative given in \equationref{secondderiv} follows the Leibnizian methodology, where the differentiation is done mechanically without considering variable dependencies.

Leibniz did, however, consider certain kinds of variables which map very directly to what we would consider as ``independent'' variables.
In the Leibniz conception, what we would consider an ``independent'' variable is \emph{a variable whose first derivative is considered constant}.
This leads to numerous simplifications of differentials because, if a differential is constant, by standard differential rules its differential is zero.
Therefore, if $x$ is the independent variable (using modern terminology) then that implies that $\dx$ is constant.
If $\dx$ is a constant (even if it is an infinitely small, unknowable constant), then that means that its differential is zero.
Therefore, $\hdiff{x}{2}$ and higher differentials of $x$ reduce to zero, simplifying the equation.\footnote{As a way of understanding this, imagine the common independent variable used in physics, especially prior to relativity---time.  
  Especially consider the way that time flows in a pre-relativistic era.  
  It flows in a continual, \emph{constant} fashion.  
  Therefore, if the flow of time (i.e., $\dt$) is constant, then by the rules of differentiation the second differential of time must be zero.  
  Thus, an independent variable is one which acts in a similar fashion to time.  
  Another way to consider this is to consider the \emph{independence} of the independent variable.  
  It's changes (i.e., differences) are, by definition, independent of anything else.   
  Therefore, we may not assign a rule about the differences between the values.  
  Thus, because there is no valid rule, the second differential may not be zero, but it is at most undefinable by definition.}

As an example, given the equation
\[
xy = 3
\]
the first differential of this would be given by
\[
x\,\dy + y\,\dx = 0
\]
and the second differential of this would be given by
\[
x\,\hdiff{y}{2} + 2\,\dx\,\dy + y\,\hdiff{x}{2} = 0.
\]
Then, you could simplify the equation by choosing any single differential to hold constant.
This is referred to in Leibnizian thought as choosing a ``progression of variables,'' and it is identical to choosing an independent variable \citep[pg.~71]{bos1974}.
Therefore, if one chooses $x$ as the independent variable, then $\dx$ is constant, and therefore $\hdiff{x}{2} = 0$. 
Thus, the equation reduces to 
\[
x\,\hdiff{y}{2} + 2\,\dx\,\dy = 0.
\]
However, if $y$ is the independent variable, then $\dy$ is held constant and therefore its differential, $\hdiff{y}{2} = 0$.
This leads to the equation
\[
2\,\dx\,\dy + y\,\hdiff{x}{2} = 0.
\]
This understanding explains the success of the modern notation of the second derivative.
The notation given in \equationref{secondderiv} is
\[
\derivop{x}{2} = \secondderiv{y}{x}.
\]
However, if we assume that $x$ is truly the independent variable, then this means that $\hdiff{x}{2} = 0$ and therefore the whole expression $\frac{\dy}{\dx}\frac{\hdiff{x}{2}}{\dx^2}$ reduces to $0$ as well.
This reduces \equationref{secondderiv} to the modern notation of $\frac{\hdiff{y}{2}}{\dx^2}$.
Additionally, if we take the assumption that $x$ is the independent variable, then the problems identified in Section~\ref{secproblem} disappear, because $x$, as an independent variable, cannot then be dependent on $t$.\footnote{To 
be clear, there is nothing preventing someone from making an independent variable dependent on a parameter.  
However, doing so then brings them around to needing to use the form of the second derivative defined here (which does not presume a particular choice of independent variable), or a compensating mechanism such as Fa\`a di Bruno's formula.
} 

In addition to \equationref{secondderiv} being reducible to $\frac{\hdiff{y}{2}}{\dx^2}$ under the assumption that $x$ is the independent variable, the Leibnizian view also gives a set of tools that allows us to reinflate instances of $\frac{\hdiff{y}{2}}{\dx^2}$ into \equationref{secondderiv}.
Euler showed that, given an equation from a specific ``progression of variables'' (i.e., a particular choice of an independent variable), we can modify that equation in order to see what it would have been if no choice of independent variable had been made.
According to \citep[pg.~75]{bos1974}, the substitution for reinflating a differential from a particular progression of variables (i.e., a particular independent variable) into one that is independent of the progression of variables (i.e., no independent variable chosen), an expansion practically identical to \equationref{secondderiv} can be used.

\section{Future Work}
The notation presented here provides for a vast improvement in the ability for higher order differentials to be manipulated algebraically.  

This improved notation yields several potential areas for study.
These include:
\begin{enumerate}
\item developing a general formula for the algebraic expansion of higher derivatives,
\item identifying additional second order differential equations that are solvable by swapping the dependent and independent variable,
\item finding other ways that differential equations can be rendered solvable using insights from the new notation, 
\item finding further reductions in special formulas that can be rendered by using algebraically manipulable notations, 
\item investigating the conjecture that the second differential of independent variables are always zero and its potential implications, and
\item extending this project to allow partial differentials to be algebraically manipulable.
\end{enumerate}

\section{Acknowledgements}

The authors would like to thank Aleks Kleyn, Chris Burba, Daniel Lichtblau, George Monta\~nez, and others who read early versions of this manuscript and provided important feedback and suggestions.

\bibliographystyle{ieeetr}
\bibliography{References}

\end{document}